\documentstyle[12pt,amsfonts]{article}
\input xypic
\title{A note on Saleh's paper `Almost continuity implies closure
continuity'\thanks{This note was written while Takashi {\sc
Noiri} was visiting the Department of Mathematics at University
of Helsinki in August 1998.\protect\newline Research supported 
partially by the Japan - Scandinavia Sasakawa Foundation.}}
\author{Julian {\sc Dontchev} and Takashi {\sc Noiri}}
\date{}
\begin{document}
\baselineskip=20pt plus 1pt minus 1pt
\newcommand{\fxy}{$f \colon (X,\tau) \rightarrow (Y,\sigma)$}
\maketitle

Recently, Saleh \cite{S1} claimed to have solved `a long standing
open question' in Topology; namely, he proved that every almost
continuous function is closure continuous (=
$\theta$-continuous). Unfortunately, this problem was settled
long
time ago and even a better result is known. Consider the
following implications:

\begin{center}
Cont. $\Rightarrow$ Almost cont. $\Rightarrow$ Almost
$\alpha$-cont. $\Rightarrow$ $\eta$-cont. $\Rightarrow$
$\theta$-cont. $\Rightarrow$ Weakly cont. \end{center}

The first two implications are trivial. In \cite{N1}, the second
author showed that every almost $\alpha$-continuous function is
$\eta$-continuous. In \cite{D1}, Dickman, Porter and Rubin proved
that every $\eta$-continuous is $\theta$-continuous and hence
weakly continuous.

Although $\theta$-continuous functions behave, in general,
nicely, they may cause some unexpected problems. For example, if
$f \colon X \rightarrow Y$ is $\theta$-continuous, then $f \colon
X \rightarrow f(X)$ is not necessarily $\theta$-continuous. Also,
the set of all points of continuity of $f \colon X \rightarrow
Y$ may be dense in $X$ and $f$ may not be weakly continuous in
either of those points. We show that with the following example:

{\em Example.} Consider the classical Dirichlet function $f
\colon ({\mathbb R},\tau_{d}) \rightarrow ({\mathbb
R},\tau_{d})$, where $\mathbb R$ is the real line with the
density topology $\tau_d$:

\[ f(x) = \left\{ \begin{array}{ll} 1, &
\mbox{$x \in {\mathbb Q},$} \\ 0, & \mbox{otherwise.}
\end{array} \right. \]

It is easily observed that $f$ is continuous at every irrational
point (furthermore, $\mathbb P$, the set of all irrationals is
dense in the density topology). On the other hand $f$ is not
weakly continuous (and hence not $\theta$-continuous) at any
irrational point.

\
Department of Mathematics, University of Helsinki, PL 4,
Yliopistonkatu 15, 00014 Helsinki, Finland, E-mail: {\tt
dontchev@cc.helsinki.fi}
\newline
Department of Mathematics, Yatsushiro College of Technology, 2627
Hirayama shinmachi, Yatsushiro-shi, Kumamoto-ken, 866 Japan,
E-mail: {\tt noiri@as.yatsushiro-nct.ac.jp}
\
\end{document}